%%%%%%%%%%%%%%%%%%%%%%%%%%%%%%%%%%%%%%%%%%%%%%%%    
%
%        THIS IS A  PLAIN TeX FILE
%
%%%%%%%%%%%%%%%%%%%%%%%%%%%%%%%%%%%%%%%%%%%%%%%%

\magnification=1200

%\font\AAA=Times at 12pt
%\font\BBB=Times at 8pt

\font\AAA=cmr14 at 12pt
\font\BBB=cmr14 at 8pt

\overfullrule=0in

\def\supp{{\rm supp}\,}
\def\Link{{\rm Link}}

\def\PSH{{\cal PSH}}
\def\dbar{\overline{\partial}}
\def\wh{\widehat}

%\def\Theorem#1{\medskip\noindent {\AAA T\BBB HEOREM \bf #1.}}
%def\Prop#1{\medskip\noindent {\bf Proposition #1.}}
%\def\Cor#1{\medskip\noindent {\bf Corollary #1.}}
%\def\Lemma#1{\medskip\noindent {\bf Lemma #1.}}
%\def\Remark#1{\medskip\noindent {\bf Remark #1.}}
%\def\Note#1{\medskip\noindent {\bf Note #1.}}
%\def\Def#1{\medskip\noindent {\bf Definition #1.}}
%\def\Claim#1{\medskip\noindent {\bf Claim #1.}}

\def\Theorem#1{\medskip\noindent {\AAA T\BBB HEOREM \rm #1.}}
\def\Prop#1{\medskip\noindent {\AAA P\BBB ROPOSITION \rm  #1.}}
\def\Cor#1{\medskip\noindent {\AAA C\BBB OROLLARY \rm #1.}}
\def\Lemma#1{\medskip\noindent {\AAA L\BBB EMMA \rm  #1.}}
\def\Remark#1{\medskip\noindent {\AAA R\BBB EMARK \rm  #1.}}
\def\Note#1{\medskip\noindent {\AAA N\BBB OTE \rm  #1.}}
\def\Def#1{\medskip\noindent {\AAA D\BBB EFINITION \rm  #1.}}

\def\Ex#1{\medskip\noindent { \AAA E\BBB XAMPLE \rm    #1.}}

\def\pf{\medskip\noindent {\bf Proof.}\ }
\def\qed{\hfill  $\vrule width5pt height5pt depth0pt$}

\def\V{W}

\def\ce{{\cal E}}

\def\vf{\varphi}

\def\and{\qquad {\rm and} \qquad}

\def\ol{\overline}
\def\bbr{{\bf R}}
\def\bbc{{\bf C}}

\def\bbz{{\bf Z}}

\def\a{\alpha}

\def\G{\Gamma}

\def\PH#1{\widehat {#1}}

\def\dbar{\ol{\partial}}
\def\wt{\widetilde}

\def\supp{\rm spt}

\  \vskip .2in
\centerline{\bf   REMARKS ON THE ALEXANDER-WERMER THEOREM  }
\smallskip

\centerline{\bf
 FOR CURVES}

\vskip .2in
\centerline{\bf F. Reese Harvey and H. Blaine Lawson, Jr.$^*$}
\vglue
2cm\smallbreak\footnote{}{ $ {} \sp{ *}{\rm Partially}$  supported by
the N.S.F. }

\centerline{\bf Abstract} \medskip
  \font\abstractfont=cmr10 at 10 pt

{{\parindent=.5in\narrower\abstractfont \noindent 
We give a new  proof of   the Alexander-Wermer  Theorem
that characterizes the oriented  curves in $\bbc^n$ which
bound positive holomorphic chains, in terms of  the linking numbers 
of the curve with algebraic cycles in the complement. 
In fact we establish  a slightly stronger version which applies
to a wider class of boundary 1-cycles.
Arguments here are based
on the Hahn-Banach Theorem and some geometric measure theory.
Several ingredients in the original proof have been eliminated.

 }}

 \vskip .5in

\centerline{\bf Table of contents.}
\vskip .3in

\moveright .6in\vbox{
1. The Alexander-Wermer Theorem 
 \smallskip
 
 2. A Dual Interpretation
 \smallskip
 
 3. The Remainder of the Proof.

}

\vfill\eject

\centerline{\bf 1. The Alexander-Wermer Theorem}
\medskip
  
 We present a different proof of the Alexander-Wermer Theorem [AW], [W$_2$]
 for curves which uses the Hahn-Banach Theorem and techniques of geometric 
 measure theory. Several ingredients of the original proof are eliminated, such 
 as the reliance on the result in [HL$_1$] that if a curve satisfies the moment
 condition, then  it bounds a holomorphic 1-chain. The arguments given here
 have been adapted to study the analogous problem in general projective manifolds
 (cf. [HL$_{3,4,5}$]).
 
 Our arguments will also apply to a more general class of curves which we now introduce.

\Def{1.1}  Let $X$ be a complex manifold and suppose  there exists a closed subset $\Sigma(\G)$ of Hausdorff 1-measure zero and an oriented, properly embedded $C^1$-submanifold of $X-\Sigma(\G)$ 
with connected components $\G_1, \G_2, ... \!$.  If, for given integers $n_1, n_2, ...$,
the sum
$\G = \sum_{k=1}^\infty n_k \G_k$ defines a current of locally finite mass in $X$
which is $d$-closed (i.e., without boundary), and if $\supp \G$ has only a finite number of connected components\footnote{$^1$}{More generally we need only assume that $\supp\G$ is contained in a compact connected set of finite linear measure} , then $\G$ will be called a 
{\bf  scarred 1-cycle (of class $C^1$) in $X$}.  By a unique choice of orientation
on each $\G_k$ we may assume each $n_k>0$.

\Ex{1.2}  Any real analytic 1-cycle is automatically a scarred 1-cycle 
(of class $C^\infty$) -- see [F, p. 433].

\Def{1.3}  By a {\bf positive holomorphic 1-chain  with boundary $\G$} we mean a sum
$V=\sum_{k=1}^\infty m_k [V_k]$  with $m_k\in \bbz^+$ and $V_k$ an irreducible
1-dimensional complex analytic subvariety of $X-\supp \G$ such that $V$ has locally finite mass in $X$ and, as currents in $X$,
$$
dV=\G %\qquad {\rm as\  currents\  in\ } X.
$$

\Remark{1.4} Standard projection techniques (cf. [Sh], [H]) show that any 1-dimensional complex  subvariety $W$ of $X-\supp \G$ automatically has locally finite 2-measure at points of $\G$,  and 
furthermore, its current boundary is of the form $dW=\sum \epsilon_k \G_k$ where
$\epsilon_k=-1, 0$ or 1 for all $k$.
See [H] and the ``added in proof'' for the more general case where $T$ is a positive $d$-closed current on  $\bbc^2-\supp\G$.

\Def{1.5}  A scarred 1-cycle $\G$ in $\bbc^n$ satisfies the  {\bf (positive) winding condition}
if
$$
{1\over 2\pi i} \int_\G {dP\over P} \ >\ 0
$$
for all polynomials $P\in \bbc[z]$ with $P\neq 0$ on $\G$.
\medskip

There are many equivalent formulations of this condition. We mention three.

\Prop{1.6}  {\sl  $\G$ satisfies the (positive) winding condition if and only if any of the
following equivalent conditions holds:

\medskip
1)\ \ 
$\int_\G d^c \varphi \ \geq\ 0$
for all smooth plurisubharmonic functions $\vf$ on $\bbc^n$.

\medskip
2)\ \ 
For each polynomial $P\in\bbc[z]$, the unique compactly supported solution $W_P(\G)$ to the equation $dW_P(\G) = P_*(\G)$ satisfies $W_P(\G)\geq 0$.

\medskip
3)\ \ The linking number $\Link(\G, Z)\geq 0$ for all algebraic hypersurfaces $Z$  contained in $ \bbc^n-\supp \G$.
}

\Prop{1.7} {\sl  If $\G$ is the boundary of a positive holomorphic 1-chain $V$ in $\bbc^n$,
then $\G$ satisfies the positive winding condition.}

\pf  We have $\int_\G d^c\vf = \int_{dV}d^c\vf= \int_{V}dd^c\vf\geq0$ since $dd^c\vf\geq0$.\qed

\medskip
The following converse of Proposition 1.7 is due to Alexander and Wermer [AW], [W$_2$].

\medskip
\noindent 
{\AAA M\BBB AIN \AAA T\BBB HEOREM} 1.8.  {\sl Let  $\G$ be a scarred 1-cycle in   $\bbc^n$.  If $\G$  satisfies the (positive) winding condition, then $\G$ bounds a positive holomorphic 1-chain in $\bbc^n$.}

\medskip
This slightly generalizes the theorem in [AW] which applies only to smooth oriented curves. However, the essential point of this paper is to provide a conceptually different proof of the result which has other applications. This proof has two distinct parts which
constitute the following two sections.

\Note{}   We adopt the following notation throughout the paper.
The polynomial hull of a compact subset $K\subset \bbc^n$ is denoted by $\PH K$.
The mass of a current $T$ with compact support in $\bbc^n$ is denoted by $M(T)$.

%%%%%%%%%%%%%%%%%%%%%%%%%%%%%%%%%%%%%%%%%%%%%%
%%%%%%%%%%%%%%%%%%%%%%%%%%%%%%%%%%%%%%%%%%%%%%
%%%%%%%%%%%%%%%%%%%%%%%%%%%%%%%%%%%%%%%%%%%%%%

\vfill\eject

\centerline{\bf 2. A Dual Interpretation}
\medskip

In this section we shall use the Hahn-Banach Theorem to establish a dual interpretation of the positive winding condition.  The main result is the following.  Recall that if $C$ is a convex cone in a topological vector space $V$, its {\bf polar} is the set
$C^0 = \{v'\in V' : v'(v)\geq0$ for all $v\in C\}$.

\Theorem{2.1. (The Duality Theorem)} {\sl The cone $A$ in the space $\ce^1_\bbr(\bbc^n)$ of smooth 1-forms on $\bbc^n$, defined by
$$
A\ \equiv \ \{d\psi+d^c \vf : \psi\in C^\infty(\bbc^n) {\rm \  and\ } \vf \in \PSH(\bbc^n)\}
$$
and the cone $B$ in the dual space $\ce_1'(\bbc^n)_\bbr$ of compactly supported one-dimensional currents in $\bbc^n$, defined by
$$
B\ \equiv \ \{S : S=d(T+R) {\rm \ with\ } T\geq 0 \ {\rm and\ } R\ {\rm of \ bidimension\ } 2,0 + 0,2\}
$$
are each the polar of the other.

Moreover, 

\smallskip

\noindent
(i)\ \ The cone $B$ coincides with the cone
$$
B'\ \equiv\ \{S : dS=0 {\rm \ and \ }\exists T\geq0 {\rm \ with\ compact\ support\  and\ } dd^cT = -d^cS\}, \
\ {\rm and}
$$
 
 \noindent
(ii)\ \ If $S=d(T+R) \in B$ with $T$ and $R$ as above, then }
$$
\supp T \ \subseteq \wh{\supp S}\ \ \ {\sl (the\ polynomial \ hull\ of\ } \supp S).
$$

This result can be restated as follows.

\Theorem{$ 2.1'$}
 {\sl A real 1-dimensional current $S$ with $dS=0$ and compact support in    $\bbc^n$ satisfies the (positive) winding condition if and only if 
$$
S\ =\ d(T+R)
\eqno{(2.1)}
$$
where $T$ is a positive 1,1 current and $R$ has bidimension  2,0 + 0,2, or equivalently,
$$
dd^c T\ =\ -d^c S
\eqno{(2.2)}
$$
for some compactly supported $T\geq 0$. 
Moreover, for each such $T$,
$$
\supp T\ \subseteq \wh{\supp S}
\eqno{(2.3)}
$$
}

\pf We will show that $B^0\ =\ A$ and that $B$ is closed.  
This is enough to conclude that $A$ and $B$ are each the polar of the other because of the bipolar theorem: $(C^0)^0 = \overline C$.

\medskip
\noindent
{\bf Proof that $A=B^0$.}  The inclusion $A\subseteq B^0$ is essentially a restatement of 
Proposition 1.7 -- the same proof applies. We need only show $B^0\subseteq A$.
Suppose $\a\in B^0$, i.e., $S(\a)\geq 0$ for all $S\in B$. Restricting to $S$ of the form $S=dR$ where $R=R^{2,0}+R^{0,2}$ is of bidimension 2,0+0,2, we see that $S(\a)=dR(\a)$ must vanish (since $-dR$ is also in $B$). Hence, $\partial\a^{1,0}=0$ and 
$\dbar \a^{0,1}=0$.  That is, $d\a=d^{1,1}\a$.  In particular, $d^{1,1}\a$ is $d$-closed.
Therefore, on $\bbc^n$ the equation $d\a = d^{1,1}\a = d d^c \vf$ can be solved for some 
$\vf \in C^\infty(\bbc^n)$.

Taking $S=dT$ where $T=\delta_p\xi\geq 0$ for $p\in\bbc^n$, yields
 $(d\a )(\delta_p\xi)=(d^{1,1}\a )(\delta_p\xi)\geq 0$.  Hence, $dd^c\vf = d^{1,1}\a\geq0$, i.e., 
$\vf\in \PSH(\bbc^n)$. Since $\a-d^c\vf$ is $d$-closed, there exists $\psi \in C^\infty(\bbc^n)$
with $\a=d\psi+d^c\vf$.  \qed
\medskip

To show that $B$ is closed requires several preliminary results.

\medskip
\noindent
{\bf Proof of (i).}
If $S\in B$, then $dd^cR$ is of bidegree $(n-1,n+1) + (n+1,n-1)$, and hence it must vanish.  Therefore, $dd^cT = -d^cS$, i.e., $S\in B'$. Conversely, if $S\in B'$, then $S-dT$ is $d^c$-closed and of course also $d$-closed. 
Note that for $T\geq 0$ and $R$ real and of bidimension 
$(2,0) + (0,2)$, the equations
$$
d(T+R)\ =\ S
\eqno{(2.4)}
$$
and
$$
\partial T + \dbar R^{n,n-2}\ =\ S^{n,n-1}
\eqno{(2.5)}
$$
are equivalent. Now the right hand side of the equation $\dbar R^{n,n-2} = 
S^{n,n-1} - \partial T $ is $\dbar$-closed.  On $\bbc^n$, this implies that there exists 
a solution $R$  with compact support.  \qed

\medskip
\noindent
{\bf Proof of (ii).}
Since $T\geq0$,   we know from [DS]  that $\supp T \subseteq
\wh{\supp\, dd^c T} $. Of course $\supp\, dd^cT=\supp\, d^cS\subseteq{\supp S}$.\qed

\Lemma {2.2}  {\sl If $S=d(T+R)\in B$, then the mass  $M(T)=T(dd^c|z|^2) = S(d^c|z|^2)$.}

\pf Note that $T(dd^c|z|^2) = (T+R)(dd^c|z|^2) = (d(T+R), d^c|z|^2)= S(d^c|z|^2)$.\qed

\Prop{2.3} {\sl The cone $B$ is closed.}
\pf
Suppose $S_j=d(T_j+R_j) \in B$ and  $S_j\to S$.  Then by Lemma 2.2, 
$M(T_j)=S_j(d^c|z|^2) \to S(d^c|z|^2)$, and so   the masses $M(T_j)$  are uniformly bounded in $j$. The convergence of $\{S_j\}$ means that all $\supp S_j\subset B(0,R)$ 
for some $R$.  Hence by part (ii) we have $\supp T_j\subset B(0,R)$ for all $j$.
By the basic compactness property of positive currents, there is a subsequence
with $T_j\to T\geq0$. Finally, since $dd^cT_j=-d^cS_j$ we have  $dd^cT=-d^cS$.
Hence, $S\in B'=B$.\qed

\vfill\eject

%%%%%%%%%%%%%%%%%%%%%%%%%%%%%%%%%%%%%%%%%%%%%
%%%%%%%%%%%%%%%%%%%%%%%%%%%%%%%%%%%%%%%%%%%%%
%%%%%%%%%%%%%%%%%%%%%%%%%%%%%%%%%%%%%%%%%%%%%

\centerline{\bf The Remainder of the Proof of the Main Theorem}
\medskip

Suppose now that $\G$ is a scarred 1-cycle in $\bbc^n$  which satisfies the 
positive winding condition.  Applying Theorem 2.1 (in its second, ``restated'' form) with $S=\G$, 
there exists  a compactly supported, positive (1,1)-current $T$ such 
$$
d(T+R)\ =\ \G
\eqno{(3.1)}
$$
where $R$ is a current of compact support and bidimension (2,0)+(0,2).
We shall show that $R=0$ and $T$ is a positive holomorphic chain.
To proceed we utilize a fundamental result of Wermer [W$_1$] in a generalized form due 
to Alexander [A].

\Theorem {3.1}  {\sl  Let $\G$ be a scarred 1-cycle of class $C^1$ in $\bbc^n$.  Then
$\wh {\supp \G} - \supp \G$  is a 1-dimensional complex analytic subvariety of  $\bbc^n-\supp \G$.}

\pf  Alexander  proves in [A] that if $K\subset\bbc^n$ is contained in a compact connected set of finite linear measure, then $\wh K-K$ is a 1-dimensional complex analytic subvariety of 
$\bbc^n-K$.  
The set $\supp \G$ has finite linear measure and only finitely many connected components. One sees from the definition that it is possible to   make a connected set $K=\supp\G\cup \tau$ of finite linear measure by adding   a   finite union of piecewise linear arcs $\tau$  contained  
  in the complement of  $\supp\G$.  
Each  irreducible component $W$  of the complex analytic curve $\wh K-K$  will have locally finite   2-measure at points of $\tau$ and will extend to $\bbc^n-\supp\G$ as a variety with  boundary 
of the form $\sum_k c_k \tau_k$, where the $c_k$'s are constants and $\tau_k$ are the connected arcs
comprising $\tau$    (cf. [HL$_1$], [H]).
Suppose this boundary is non-zero.
Then $W$ must be contained in the union of the complex lines determined by the real line
segments comprising $\partial W \cap \tau$. Since $W$ is irreducible,  it is contained in just
one such complex line. Constructing $\tau$ so that each connected component of $\tau$ has at least
two (complex independent) line segments, we have a contradiction. Thus, for generic choice
of $\tau$, the set   $\PH K-\supp\G$ is a 1-dimensional subvariety of $\bbc^n-\supp\G$. 
 In particular, this proves that   $\PH K\subseteq \PH {\supp\G}$. 
 Since $\PH {\supp\G}\subseteq \PH K$, we are done.
\qed\medskip

Let $V_1,V_2,... $ denote the irreducible components of the complex curve given by Theorem 3.1.  
We are going to prove that $T=\sum_j n_j V_j$ for positive integers  $n_j$.  For this  we first utilize a result from [HL$_2$, p. 182].

\Lemma{3.2}  {\sl Suppose $T$ is a positive current of bidimension 1,1 with $dd^c T=0$ on a complex manifold $X$.  If $T$ is supported in a complex analytic curve $W$ in $X$, then $T$ can be written as a sum $T=\sum_j h_j W_j$ where each $W_j$ is an irreducible component of $W$ and $h_j$ is a non-negative harmonic function on $W_j$.}
\medskip
The case needed here is the following.

\Cor{3.3}  {\sl  If $T\geq 0$ satisfies $dd^cT = -d^c\G$ on $\bbc^n$, then on $\bbc^n-\supp\G$
one has $T=\sum_j h_j V_j$ with $h_j$ harmonic on $V_j$.}
\medskip

We first restrict attention to dimension $n=2$, where the equation (2.5), namely
$$
\dbar R^{2,0}\ =\ \G^{2,1} - \partial T
$$
implies that $R^{2,0}$ is a holomorphic 2-form outside the support of $\G-dT$.

\Lemma{3.4. (n=2)}  {\sl  If $d(T+R)=\G$ with $T\geq 0$ and $R$ of bidimension
$(2,0) + (0,2)$, then}
$$
\supp R\ \subseteq \ \wh{\supp \G}
$$

\pf
By Theorem 2.1(ii),  $R^{2,0}$ is a holomorphic 2-form on $\bbc^2-\wh{\supp\G}$, and $R^{2,0}$ vanishes outside
of a compact subset of $\bbc^2$.  The polynomially convex set $\wh{\supp\G}$ cannot have a bounded component in its complement.  Therefore, $R^{2,0}$  must vanish on all of $\bbc^2-\wh{\supp\G}$.
\qed

\Lemma{3.5. (n=2)}  {\sl Each $h_j\equiv c_j$ is constant, and the current $T=\sum_j c_j V_j$ is $d$-closed on $\bbc^2-\supp\G$.}

\pf  Pick a regular point of one of the components $V_j$, let $\pi$ denote a holomorphic projection
(locally near the point) onto $V_j$, and let $i$ denote the inclusion of $V_j$ into $\bbc^2$.  
Note that $T$ is  locally supported in $V_j$ by Theorem 2.1(ii) while $R$ is locally supported in $V_j$ by Lemma 3.4. Therefore, both of the push-forwards $\pi_*T$ and $\pi_*R$ are well defined.  Now
$\pi_*R$, being of bidimension $(2,0) + (0,2)$ in $V_j$ must vanish.  However, $T=h_jV_j$ satisfies
$\pi_*T = h_j$.  Since $d(T+R)=0$, the push-forward $\pi_* d(T+R)= d\pi_* (T+R)= dh_j$ must also vanish, i.e., each $h_j = c_j$ is constant.  This proves:

\def\V{T^0}

\Cor{3.6. (n=2)}  {\sl  The current $T=\sum_j c_j V_j$ on $\bbc^2-\supp\G$ has locally finite mass across
$\supp\G$ and its extension $\V$ by zero across $\supp\G$ satisfies
$$
d\V\ =\ \sum_j r_j  \G_j  \qquad {\rm \ on\ \ } \bbc^2
$$
for real constants $r_j$.}
\pf   See Remark  1.4.\qed\medskip

Another corollary of Lemma 3.5 is the following.

\Cor{3.7}  \qquad $\supp R\ \subseteq\ \supp \G$

\pf   By (2.5) the current $R^{2,0}$ is a holomorphic 2-form on $\bbc^2-\supp\G$
since $dT=0$ there. Since $R^{2,0}$ vanishes at infinity, this proves the result. \qed\medskip

\noindent
{\bf Completion of the case n=2.}  Now
$$
T+R \ =\ \V+\chi T+R
\eqno{(3.2)}
$$
where $\chi$ is the characteristic function of $\supp\G$ and $\chi T +R$ has support in $\supp\G$.
We also have
$$
d(\chi T+R)\ =\ \sum_j(n_j-r_j)\G_j\qquad{\ \rm on\ \ }\bbc^2
\eqno{(3.3)}
$$
Let $\rho$ denote a local projection onto a regular point of $\G_j$.
Then $\rho_*(\chi T+R)$ is a well defined current on $\G_j$, but of dimension 2.
Hence it must vanish. 
Since $\rho_*$ commutes with $d$, this proves that $(n_j-r_j)\G_j$ must vanish.
Hence, $r_j=n_j$ for all $j$,  and so  $d(\chi T+R)=0$ and $d\V = d(T+R)$ by equations (3.2) and 
(3.3).
 This proves that  $d\V=\G$  by (3.1) .

\medskip
\noindent
{\bf Proof for the case $ n\geq 3$.}
The general case follow easily from the case where $n=2$.  Consider a generic linear projection
 $\pi:\bbc^n\to\bbc^2$ so that each mapping $V_j\to \pi V_j$ is one-to-one.
 Then the current $T=\sum_j h_j V_j$ in $\bbc^n-\supp\G$ projects 
 to the current $\pi_*T = \sum_j \wt{h}_j \pi(V_j)$ in $\bbc^2-\pi(\supp\G)$ where $\wt{h}_j\circ \pi = h_j$.
Since  each $\wt{h}_j=c_j$ is constant, so is each $h_j$. Now the current $\V=\sum_jc_jV_j$ satisfies
$d\V=\sum_jr_j\G_j$ and again by projecting we conclude that $r_j=n_j$.\qed

\vskip .4in

% \magnification=1200
%\NoBlackBoxes
%\nologo 

% \input qtmacros 
% \input QASdefs.tex 

\centerline{\bf References}

\vskip .2in

\noindent
\item{[A]}   H. Alexander, {\sl Polynomial approximation and hulls in sets of finite linear measure in $\bbc^n$},
 Amer. J. Math. {\bf 93} (1971),  65-74.

\smallskip

% \noindent
%[AW$_1$]  H. Alexander and J. Wermer,  Several Complex
%Variables and Banach Algebras, Springer-Verlag,  New York, 1998. 

 %\smallskip

\noindent
\item{[AW]}  H. Alexander and J. Wermer, {\sl Linking numbers
and boundaries of varieties}, Ann. of Math.
{\bf 151} (2000),   125-150.

 \smallskip

 \noindent
\item{[DS]} J. Duval and N. Sibony, {\sl
Polynomial convexity, rational convexity and currents},
  Duke Math. J. {\bf 79}  (1995),     487-513.

 \smallskip

\noindent
\item{[F]}   H. Federer, Geometric Measure  Theory,
 Springer--Verlag, New York, 1969.

 \smallskip

\noindent
\item{[H]}  F.R. Harvey,
Holomorphic chains and their boundaries, pp. 309-382 in ``Several Complex
Variables, Proc. of Symposia in Pure Mathematics XXX Part 1'', 
A.M.S., Providence, RI, 1977.

\noindent
\item{[HL$_1$]}  F. R. Harvey and H. B. Lawson, Jr, {\sl On boundaries of complex
analytic varieties, I}, Annals of Mathematics {\bf 102} (1975),  223-290.

 \smallskip

\noindent
\item{[HL$_2$]} F. R. Harvey and H. B. Lawson, Jr, {\sl An intrinsic
characterization of K\"ahler manifolds}, Inventiones Math.,  {\bf 74} 
(1983), 169-198.
 \smallskip

 \noindent
\item{[HL$_3$]}  F. R. Harvey and H. B. Lawson, Jr, {\sl Projective hulls and
the projective Gelfand transformation},  Asian J. Math. {\bf 10}, no. 2 (2006), 279-318. ArXiv:math.CV/0510286.

 \smallskip

\noindent
\item{[HL$_4$]}  F. R. Harvey and H. B. Lawson, Jr, {\sl Projective linking and boundaries of positive holomorphic chains in projective manifolds, Part I},  Stony Brook Preprint, 2004.

ArXiv:math.CV/0512379

 \smallskip

 \noindent
\item{[HL$_5$]}  F. R. Harvey and H. B. Lawson, Jr, {\sl Boundaries of positive holomorphic chains},   
 Stony Brook Preprint, 2006.

 \smallskip

   \noindent
\item{[Sh]}    B. Shiffman,    {\sl  On the removal of singularities of analytic sets},    Michigan  Math. J.,
 {\bf  15}  (1968), 111-120  .

\smallskip

%   \noindent
%\item{[S] }  H. H. Schaefer,  Topological Vector Spaces,    Springer Verlag,
%New York,  1999.

%\smallskip

 \noindent
\item{[W$_1$]}   J. Wermer  {\sl    The hull of a curve in $\bbc^n$},    
Ann. of Math., {\bf  68}  (1958), 550-561.

\smallskip

 \noindent
\item{[W$_2$]}   J. Wermer    {\sl    The argument principle and
boundaries of analytic varieties},     Operator Theory: Advances and 
Applications, {\bf  127}  (2001), 639-659.

\smallskip

\end